\documentclass{article}

\usepackage{arxiv}

\usepackage[utf8]{inputenc} 
\usepackage[T1]{fontenc}    
\usepackage{hyperref}       
\usepackage{url}            
\usepackage{booktabs}       
\usepackage{amsfonts}       
\usepackage{nicefrac}       
\usepackage{microtype}      
\usepackage{lipsum}		
\usepackage{graphicx}
\usepackage[numbers]{natbib}
\usepackage{doi}

\usepackage{amsmath,amssymb,amsfonts}%
\usepackage{amsthm}%
\usepackage{mathrsfs}%
\usepackage{float}

\usepackage{multirow}
\usepackage{bm}

\usepackage{algorithm}%
\usepackage{algorithmicx}%
\usepackage{algpseudocode}%
\usepackage{listings}%

\theoremstyle{thmstyleone}%
\theoremstyle{thmstyletwo}%

\theoremstyle{thmstylethree}%

\title{Compactly Supported Radial Basis Functions as Probability Density Functions}

\author{ \href{https://orcid.org/0009-0004-9623-1766}{\includegraphics[scale=0.06]{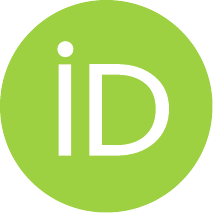}\hspace{1mm}Sergio~Díaz-Elbal} \\\\
	Department of Mathematics and\\
	Center for Development and\\
	Transfer of Mathematical Research\\
	to Industry (CDTIME),\\
	University and Almería,\\
	Almería, Spain \\\\
	\texttt{sdiazelbal@ual.es} \\
	\And
	\href{https://orcid.org/0000-0001-9421-5624}{\includegraphics[scale=0.06]{orcid.pdf}\hspace{1mm}Andrei~Martínez-Finkelshtein} \\\\
	Department of Mathematics and\\
	Center for Development and\\
	Transfer of Mathematical Research\\
	to Industry (CDTIME),\\
	University and Almería,\\
	Almería, Spain \\\\
	Department of Mathematics,\\
	Baylor University,\\    
    Waco, Texas, United States \\\\
	\texttt{andrei@ual.es} \\
	\And
	\href{https://orcid.org/0000-0002-6127-6559}{\includegraphics[scale=0.06]{orcid.pdf}\hspace{1mm}Darío ~Ramos-López} \thanks{Corresponding author.} \\\\
	Department of Mathematics and\\
	Center for Development and\\
	Transfer of Mathematical Research\\
	to Industry (CDTIME),\\
	University and Almería,\\
	Almería, Spain \\\\
	\texttt{dramoslopez@ual.es} \\    
}

\renewcommand{\shorttitle}{Compactly Supported Radial Basis Functions as Probability Density Functions}

\hypersetup{
pdftitle={Compactly Supported Radial Basis Functions as Probability Density Functions},
pdfsubject={math.ST},
pdfauthor={Sergío~Díaz-Elbal, Andrei~Martínez-Finkelshtein,Darío~Ramos-López},
pdfkeywords={Radial Basis Functions, Probability Density Function, Density Estimation, Wendland Function},
}

\begin{document}
\maketitle

\begin{abstract}
	Compactly Supported Radial Basis Functions (CS-RBFs) are a fundamental tool in multivariate approximation theory. However, their use in statistics and probability modeling remains underexplored, having been used mainly to express covariance functions in Gaussian processes or as kernel functions. This work explores CS-RBFs as a novel parametric family of probability density functions, focusing in particular on Wendland $\mathscr{C}^2$ kernels. The primary contribution of this work is the derivation of analytical expressions for various statistical properties, such as moments and the cumulative distribution function, of CS-RBFs as univariate and conditional densities. The approach comprises two alternative scenarios: when the CS-RBF support lies entirely within the variable's domain (untruncated support) and when part of it is outside (truncated support). Mixture models employing CS-RBFs are also analyzed, and their main properties are detailed. Furthermore, we introduce an incremental learning algorithm for density estimation with CS-RBF mixture models, in which centers are determined using k-means and weights and shape parameters are optimized by stochastic gradient descent. Experiments on synthetic and real-world datasets show that CS-RBF densities provide competitive results in terms of likelihood and model complexity in comparison with Gaussian mixture models. In addition, these CS-RBF densities allow the exact computation of key distributional properties in univariate and conditional settings.
\end{abstract}


\section{Introduction}\label{intro}

Density estimation \citep{Silverman1986,Scott2015} is a fundamental task in statistics and machine learning, serving as the basis for numerous applications across various domains, such as anomaly detection \citep{Liu2022}, clustering \citep{Scrucca2023} and imaging \citep{Rajwade2008}. Its primary goal is to reconstruct an unknown probability density function from a finite set of independent and identically distributed (i.i.d.) samples, providing a faithful model of the underlying uncertainty of the continuous random variable.

Nonparametric methods, such as kernel density estimation (KDE), are widely employed due to their simplicity and effectiveness, particularly when no prior knowledge of the distributional form is available \citep{Cattaneo2020,Cattaneo2022,Gramacki2018}. KDE can adapt to several different geometries by placing a kernel function at every data point. Since the KDE model complexity scales linearly with the sample size, it can become computationally demanding for large datasets in terms of both memory storage and computing time \citep{Gramacki2018,Charikar2023}.

Furthermore, in high-dimensional spaces, standard KDE techniques suffer from the curse of dimensionality \citep[Chap.~3]{Gramacki2018}, often requiring large amounts of data to produce reliable estimates. Moreover, kernel methods often suffer from boundary bias; for instance, when the data are non-negative, they assign probability mass outside the support of the variable, since these methods do not incorporate prior information about the data boundaries \citep{Bertin2019,Cattaneo2020}. This problem often requires the use of complementary techniques \citep{Colbrook2020}, such as boundary reflection methods.

Finite mixture models \citep{McLachan2019}, such as Gaussian mixture models (GMMs) \citep{Rasmussen1999,Scrucca2023}, offer a better compromise between flexibility and complexity, when compared to nonparametric methods. By representing the density as a weighted sum of a fixed number of component densities, mixture models constitute a semi-parametric approach to density estimation. They can be seen as parametric, since each component density depends on some parameters, but they can be seen as nonparametric because the number of components can be chosen arbitrarily. GMMs, in particular, benefit from the universality of the Gaussian distribution and well-established fitting algorithms \citep[Chap.~6]{Bishop2006}, such as Expectation-Maximization (EM). However, a significant limitation of standard GMMs is the infinite support of the Gaussian basis. When applied to bounded continuous variables, Gaussian components inevitably assign non-zero probability mass outside of the domain. This phenomenon, known as boundary bias, occurs because the model does not account for the natural bounds of the data. While Truncated Gaussian Mixture Models \citep{Lee2012} can address this by confining the support, they may incur additional computational overhead because they require normalization constants, which usually must be computed numerically. Furthermore, GMMs are highly sensitive to outliers \citep{Kenig2024}; because the model relies on mean and covariance, extreme values can inflate the variance of the components, potentially degrading the quality of fit. 

Compactly Supported Radial Basis Functions (CS-RBFs) \citep{Buhmann2001,Wu1995,Wendland2005} have become a cornerstone in approximation theory and numerical analysis, due to their diverse applications, such as meshfree methods for PDEs \citep{Wong2002}, interpolation and approximation \citep{Fasshauer2007}, surface reconstruction \citep{Morel2018}, medical imaging \citep{Mohammadi2019} or dynamical models in physics \citep{Baes2024, Dehnen2012}.

In this work, we aim to connect approximation theory and probabilistic modeling by exploring the use of CS-RBFs as a parametric family of probability density functions, and derive their main distributional properties. The use of CS-RBFs in the context of probability and statistics is rather limited. Some works in this topic employ CS-RBFs as covariance kernels in Gaussian Process Regression \citep{Risser2025} or as smoothing kernels in the fashion of KDE estimation \citep{Dehnen2012}. An exception is \citet{Lin2006},An exception is Lin and Yuan (2006), who use CS-RBFs as reproducing kernels for
regularization in a white noise model, a setting asymptotically equivalent to density estimation, rather than as densities in their own right. Therefore, the use of CS-RBFs to model probability distributions or to serve as a basis for more complex statistical models remains largely unexplored.

Thus, our first goal is to analyze the theoretical properties of CS-RBFs as a family of parametric probability density functions. We consider univariate, conditional and multivariate distributions. Our approach covers two alternative cases: when the CS-RBF support lies entirely within the variable's domain (untruncated support) and when part of it is outside (truncated support). This study also addresses mixture models with CS-RBFs. Although the focus is on the Wendland $\mathscr{C}^2$ function, the methodology and most derived properties can be extended to other CS-RBFs families.
 
The structure of the paper is as follows. In Section \ref{sec:csrbfs}, we review our basic tool for building CS-RBFs: the Wendland $\mathscr{C}^2$ function. We explore how it is constructed with truncated support, to ensure no probability is assigned outside of a continuous variable domain. To provide a rigorous theoretical foundation to CS-RBF densities, Sections \ref{sec:unipdf} and \ref{sec:condpdf} detail the derivation of their key analytical properties, including closed-form expressions for moments, normalization constants, and the cumulative distribution function (CDF) for both univariate and conditional densities with truncated or untruncated domains. Based on these properties, Section \ref{sec:mixtures} introduces CS-RBF mixture models. As an application, in Section~\ref{csrbfs:densityestimation} we propose a learning procedure for a CS-RBF mixture density. This approach consists of an incremental algorithm that incorporates Stochastic Gradient Descent (SGD) to obtain the mixture parameters. In it, we employ reparameterization techniques to ensure that all parameters remain within their appropriate ranges. Finally, in Section \ref{sec:experiments}, we validate the learning algorithm for CS-RBF mixture density estimation experimentally by fitting synthetic and real-world data in the univariate, conditional, and multivariate cases.

\section{Compactly Supported Radial Basis Functions} \label{sec:csrbfs}

Let $\mathbb{R}^d$ denote the $d$-dimensional Euclidean space equipped with the standard Euclidean norm $\|\cdot\|$. A multivariate function $\Phi: \mathbb{R}^d \rightarrow \mathbb{R}$ is called radial if there exists a univariate function $\phi: [0, \infty) \rightarrow \mathbb{R}$ such that for all $x \in \mathbb{R}^d$,
\begin{equation*}
    \Phi(x) = \phi(r), \quad \text{where } r = \|x\|.
\end{equation*}
Radial Basis Functions (RBFs) are central to multivariate approximation theory, particularly for the interpolation and approximation of scattered data \citep{Wendland2005,Fasshauer2007}. While globally supported RBFs (e.g., the Gaussian $\phi(r) = e^{-r^2}$) are theoretically robust and infinitely smooth, they take positive values everywhere and have an unbounded support. In contrast, Compactly Supported Radial Basis Functions (CS-RBFs) are characterized by a bounded support, in which the function takes positive values. Consequently, $\Phi(x) = 0$ when $\|x\|$ exceeds a specific threshold, leading to sparse representations and localizing the influence of each basis function.

Prominent CS-RBF families include those constructed by  \citet{Buhmann2001}, \citet{Wu1995}, and \citet{Wendland2005}.

The Wendland family constitutes a class of CS-RBFs that are piecewise polynomials with a prescribed order of smoothness and guaranteed positive definiteness up to a maximum spatial dimension $d$. In this work, we specifically employ the Wendland $\mathscr{C}^2$ function \citep[Chap.~9]{Wendland2005}. In the literature, it is denoted by $\phi_{3,1}$, but we refer to it simply as $\phi$. This function is twice continuously differentiable and takes the form of a univariate polynomial on its support domain $[0,1]$. It is defined explicitly as follows:
\begin{equation} \label{eq:wendland_c2}
    \phi(r) = (1 - r)_+^4 (4r + 1),
\end{equation}
where $(\cdot)_+$ denotes the positive part of the function, so that $(1-r)_+ = \max(0, 1-r)$. By construction, this specific function is strictly positive definite for dimensions $d \leq 3$ \citep[Chap.~9]{Wendland2005}. Nevertheless, it is important to clarify the role of this property in our context. Positive definiteness is a standard requirement for the  Gram matrices in interpolation and approximation problems. In contrast, our framework relies solely on the function's non-negativity and integrability to form valid probability density functions. Therefore, positive definiteness is not a property we need in our statistical setting. 

The function defined in \eqref{eq:wendland_c2} provides several critical advantages for density estimation. First, because $\phi$ is a polynomial on its support, it avoids the evaluation of computationally expensive transcendental functions (such as the exponential function required by Gaussians). Second, its $\mathscr{C}^2$ continuity ensures that the density is smooth and has continuous gradients, which is required by the gradient-based optimization methods we employ for density estimation.

To adapt the support of the basis function to the local geometry of the data, we introduce a shape parameter $\varepsilon > 0$, which controls the compact support radius (the support radius being $1/\varepsilon$). Moreover, the functions can be centered at any point, which we call the center $c$. The centered and scaled CS-RBF is defined as:
\begin{equation}\label{eq:csrbf}
    \Phi_{\varepsilon,c}(x) = \phi\left(\varepsilon\|x-c\| \right).
\end{equation}
Under this transformation, the function is strictly zero for all $\|x-c\| \geq 1/\varepsilon$. This spatial scalability, combined with compact support and computational efficiency, makes the Wendland $\mathscr{C}^2$ function an ideal candidate to act as a probability density function.

\section{CS-RBFs as Univariate Probability Density Functions} \label{sec:unipdf}
 This section formulates compactly supported radial basis functions as a parametric family of probability density functions. We aim to study the main properties of these CS-RBFs as density functions. In particular, we focus on the Wendland $\mathscr{C}^2$ function because of its simplicity and continuity. However, many of the properties we derive can be computed for other functions belonging to these Wendland families. A univariate probability density function in $\Omega=[a,b]$ is given by
\begin{equation}
    \label{eq:univar_density_CSRBF}
     f(x)= \frac{1}{N} \phi\left(\varepsilon|x-c| \right),
 \end{equation}
where $N$ denotes the normalization constant required for $f$ to integrate 1 in $\Omega$. A primary motivation for selecting the Wendland $\mathscr{C}^2$ function $\phi$ is its smoothness, compact support and low-degree polynomial form. 
 
Unlike Gaussian kernels, which require error functions for integrals over finite domains \citep{Jawitz2004}, the polynomial nature of Wendland functions allows the derivation of closed-form expressions for normalization constants, moments, and cumulative distribution functions on arbitrary intervals. 

Therefore, we assume our density function $f(x)$ in \eqref{eq:univar_density_CSRBF} is defined on an arbitrary interval $\Omega = [a,b]$ and is identically zero elsewhere. Formally, it can be written as $f(x) = \frac{1}{N}\mathbf{1}_{\Omega}(x)\,\phi(\varepsilon|x-c|)$, but for clarity we can omit the indicator function $\mathbf{1}_{\Omega}(x)$ in the expression as $f(x)$ is only used within $\Omega$.

If the function's support $D=[c-\frac{1}{\varepsilon},c+\frac{1}{\varepsilon}]$ is not contained in $\Omega$, then the CS-RBF support is truncated. Thus, we distinguish between two cases: the untruncated support, when $D$ is fully contained within $\Omega$ and the truncated support, when it is cut off by $\Omega$.

\subsection{Normalization constant}

Consider the univariate CS-RBF density function defined as in \eqref{eq:univar_density_CSRBF}, with center $c$ and shape parameter $\varepsilon$. In order to obtain a proper density function, we need to be able to find an expression for the normalization $N$. That is, we must compute the exact integral of the CS-RBF in the interval $\Omega$,
\begin{equation*}
N =  \int_a^b \phi\bigl(\varepsilon|x-c|\bigr)\,dx.
\end{equation*} 
Notice that $N>0$ only when $\left( D\cap\Omega \right)^\circ \neq \emptyset$; otherwise, considering \eqref{eq:univar_density_CSRBF} as a density function in $\Omega$ is meaningless.

To that end, we define the following primitive function, which is an explicit polynomial that remains constant for $r \geq 1$ (since $\phi(r)=0$ for $r\geq1$), 
\begin{equation}
  \label{eq:F0}
  F_0(r)= \begin{cases} 
        \int_0^r \phi(s)ds = \frac{2 r^6}{3}-3 r^5+5r^4-\frac{10
   r^3}{3}+r,\quad &0\leq r \leq 1,\\
       \int_0^1 \phi(s)ds = \frac{1}{3} \quad &r>1.
    \end{cases}
\end{equation}

Given fixed $\varepsilon$ and center $c$, we define the following function
\begin{equation}\label{eq:X(x)}
    X(x):= \varepsilon|x-c|.
\end{equation}
From this function, we can define the new transformed boundaries
\begin{equation}\label{eq:limits}
    A=X(a), \quad B=X(b).
\end{equation}
Note that $A<B$ if $c<a$, and $A>B$ if $c>b$. Therefore, we can compute the normalization constant as: 
\begin{equation}\label{eq:normalization1d}
    N=\begin{cases}
    \dfrac{1}{\varepsilon} \, \left|F_0(A)-F_0(B) \right|, & c\notin [a,b], \\[3mm]
    \dfrac{1}{\varepsilon} \, \left(F_0(A)+F_0(B) \right), & c\in [a,b].
\end{cases}
\end{equation}

The case $c\notin [a,b]$ handles $c<a$ and $c>b$ through the absolute value. The case $c\in [a,b]$ is obtained by splitting the integral into two parts from the center of the function. When the support is untruncated, $\varepsilon(c-a)\geq 1$ and $\varepsilon(b-c)\geq 1$, so we can obtain an explicit normalization 
\begin{align}\label{eq:norm_full}
N = \dfrac{1}{\varepsilon} \, \left(F_0(A)+F_0(B) \right)=\dfrac{1}{\varepsilon} \, \left(F_0(1)+F_0(1) \right) =\frac{2}{3 \varepsilon}.
\end{align}

\subsection{Moments}

Among the most useful and relevant summary quantities associated with a probability density are its moments. These moments characterize a probability distribution and constitute the basis for different data analysis and fitting techniques. The availability of closed-form expressions for these moments enables their exact computation without the cumulative errors or computational costs associated with numerical quadrature techniques.

Given a density function defined as in \eqref{eq:univar_density_CSRBF}, its $k$-th moment is defined as
\begin{equation}\label{eq:moment}
m_{k} = \frac{1}{N} \int_a^b x^k\, \phi\bigl(\varepsilon|x-c|\bigr)\,dx,
\end{equation}
where $N$ is the normalization constant obtained in \eqref{eq:normalization1d}. Again, we consider the domain $\Omega=[a,b]$. Let us define the generalized primitive function
\begin{equation}\label{eq:Fk}
    F_k(r) = \begin{cases}
        \int_0^r s^k \phi(s)ds=r^{k+1} \left(\frac{4 r^5}{k+6}-\frac{15r^4}{k+5}+\frac{20 r^3}{k+4}-\frac{10
   r^2}{k+3}+\frac{1}{k+1}\right), &0\leq r\leq 1, \\
   \int_0^1 s^k(1-s)^4(4s+1)ds = \frac{120(k+2)k!}{(k+6)!} , \quad &r>1.
       \end{cases}
\end{equation}
Note that the case $k=0$ recovers the previous primitive $F_0$ in \eqref{eq:F0}.

Proceeding similarly to the computation of the normalization constant, we perform the change of variables $u = x - c$, and expand $(u+c)^k$, obtaining the following relation:
\[
m_{k} = \frac{1}{N} \int_{a-c}^{b-c} (u+c)^k\, \phi\bigl(\varepsilon|u|\bigr)\,du= \frac{1}{N} \sum_{j=0}^{k} \binom{k}{j}\, c^{\,k-j}  M_{j} \;,
\]
where $M_{j}= \int_{a-c}^{b-c} u^j\, \phi(\varepsilon  |u|)\,du$. To compute $M_{j}$, we employ the generalized primitive $F_k$ defined in \eqref{eq:Fk}. Recalling our transformed boundaries defined in \eqref{eq:limits}, $A=\varepsilon |a-c|$ and $B=\varepsilon |b-c|$, we apply the substitution $t = \varepsilon |u|$ to the integral $M_{j}$ so that it matches the primitive $F_k$. Since the integrand $u^j \phi(\varepsilon |u|)$ is even or odd depending on the parity of $j$, a factor $(-1)^j$ is needed in some cases to account for the contribution from the negative part of the domain. The shifted moment is then obtained by evaluating the corresponding $F_j$ at the transformed boundaries, accounting for whether the center lies within $[a,b]$ or not:
\begin{equation}
    M_j = \frac{1}{\varepsilon^{j+1}}
    \begin{cases}
        F_j(B) - F_j(A),                   & c < a, \\[6pt]
        (-1)^j F_j(A) + F_j(B),            & c \in [a, b], \\[6pt]
        (-1)^j \bigl(F_j(A) - F_j(B)\bigr),& c > b.
    \end{cases}
    \label{eq:centered_moments}
\end{equation}

\subsection{Cumulative Distribution Function}

Given a domain $\Omega=[a,b]$, the cumulative distribution function (CDF) of a continuous random variable with density function $f$ is defined as 

\begin{equation*}
    F(x)=\int_{-\infty}^x f(t)\,dt.
\end{equation*}

For a specific univariate CS-RBF with domain $\Omega=[a,b]$, we must compute:
\begin{equation}\label{eq:cdf_integral}
F(x) = \frac{1}{N} \int_{a}^x \phi\bigl(\varepsilon|t-c|\bigr)\,dt.
\end{equation}

Recalling our transformed boundary defined in \eqref{eq:limits}, $A = \varepsilon|a-c|$, and the function $X(x)$ defined in \eqref{eq:X(x)}, the exact analytical solution to \eqref{eq:cdf_integral} is determined directly by the geometric position of the center $c$ relative to the domain $[a,b]$. If $c \in [a,b]$ then
\begin{equation}\label{eq:cdf_global} 
F(x) = \begin{cases}  
0, & x < a, \\[3mm]
\frac{1}{N\varepsilon} \Bigl(F_0(A) + (-1)^{\delta} F_0\bigl(X(x)\bigr) \Bigr),  \quad  & a \leq x \leq b, 
\\[3mm] 1, & x > b.
\end{cases}
\end{equation}
where $\delta=0$ if $x>c$ and $\delta=1$ if $x\leq c$. If, on the contrary, $c\notin [a,b]$, then the CDF can be written as
\begin{equation}\label{eq:cdf_outside}
F(x) = \begin{cases}  
0, & x < a, \\[3mm]
\frac{1}{N\varepsilon} \Bigl|F_0\bigl(X(x)\bigr) - F_0(A)\Bigr|, & a \leq x \leq b, \\[3mm]
1, & x > b.
\end{cases}
\end{equation}

The continuity of these two functions is straightforward to verify. This formulation allows for the exact evaluation of probabilities and quantiles without numerical integration, even when the variable domain truncates the natural support of the CS-RBF density function.

\section{CS-RBFs as Conditional Density Functions} \label{sec:condpdf}

CS-RBFs are naturally defined in multivariate spaces; however, their compact support makes their application as conditional densities particularly powerful. In this section, we transition from a multivariate joint density to a univariate conditional density, a process that involves ``slicing'' the radial function and re-normalizing the resulting univariate function. From Equation \eqref{eq:csrbf} we define a CS-RBF multivariate probability density function as follows: 
\begin{equation*}
    f(\mathbf{x})=\frac{1}{N}\phi(\varepsilon\|\mathbf{x}-\mathbf{c}\|)
\end{equation*}
where $\mathbf{x},\mathbf{c} \in \mathbb{R}^d$, $\varepsilon$ is the shape parameter, $N$ denotes the corresponding normalization constant.

\subsection{Conditioning}

Let $\mathbf{X}=(X_1,\mathbf{X}_{2:d})\in \mathbb{R}^d$ denote a $d$-dimensional random vector, where $X_1$ is the target variable and $\mathbf{X}_{2:d}=(X_2,\dotsc,X_d)$ are the conditioning variables. The conditional density of $X_1$ given the values $\mathbf{X}_{2:d}=\mathbf{x}_{2:d}$ is, by Bayes' theorem, 
\begin{equation}\label{eq:bayes}
    f_{X_1|\mathbf{X}_{2:d}}(x_1| \mathbf{x}_{2:d}) = \frac{ f_{\mathbf{X}}(\mathbf{x})  }{ f_{\mathbf{X}_{2:d}}(\mathbf{x}_{2:d}) }.
\end{equation}
For brevity we write $f(x_1|\mathbf{x}_{2:d})$,
 $f(\mathbf{x})$ and $f(\mathbf{x}_{2:d})$ for the conditional, joint and marginal densities respectively.

Along the rest of Section~\ref{sec:condpdf}, we focus on the conditional density $f(x_1|\mathbf{x}_{2:d})$ of $x_1 \in \Omega = [a, b]$ given a fixed vector of conditioning values $\mathbf{x}_{2:d} = (x_2, \ldots, x_d)$, with the aim of analyzing its statistical properties. 

Let us define the Euclidean distance from the vector of conditioning values $(x_2, \ldots, x_d)$ to the corresponding center coordinates $(c_2, \hdots, c_d)$ as
\begin{equation}\label{eq:rho}
    \rho = \left(\sum_{i=2}^{d}(x_i - c_i)^2\right)^{1/2},
\end{equation} 
and the effective conditional radius as $\Delta = \sqrt{\varepsilon^{-2} - \rho^2}$.
If $\rho \geq 1/\varepsilon$, the conditioning values lie outside the support and the conditional density is identically zero. Otherwise, if $\rho<1/\varepsilon$, then $\Delta >0$ and the one-dimensional support for $x_1$ is $[c_1 - \Delta,\, c_1 + \Delta]$, which, intersected
with $\Omega$, gives the effective domain bounds
\begin{equation}
    L = \max(a,\, c_1 - \Delta),  \quad U = \min(b,\, c_1 + \Delta).
    \label{eq:LU}
\end{equation}
If $[c_1-\Delta,\,c_1+\Delta] \cap \Omega = \emptyset$, then $L > U$, and the conditional density is identically zero, as in the case $\rho \geq 1/\varepsilon$ above. The extreme case $L = U$ occurs only when this intersection is a single point, in which case the domain of $X_1 \mid \mathbf{x}_{2:d}$ collapses to a point and the conditional variable is degenerate. We henceforth exclude this case and consider $L < U$.

Applying Bayes' formula \eqref{eq:bayes} to the multivariate CS-RBF density, and decomposing the Euclidean norm as $\|\mathbf{x}-\mathbf{c}\|^2 = (x_1-c_1)^2 + \rho^2$, the slice is a one-dimensional function of $x_1$ alone. The global normalization $N$ cancels and the conditional density is simply the unnormalized one-dimensional slice renormalized over $[L, U]$:
\begin{equation}
    f(x_1 \mid \mathbf{x}_{2:d}) 
    =\frac{f(\mathbf{x})}{f(\mathbf{x}_{2:d})}= \frac{\phi\!\left(\varepsilon\sqrt{(x_1-c_1)^2+\rho^2}\right)}
           {\displaystyle\int_L^U \phi\!\left(\varepsilon\sqrt{(x-c_1)^2+\rho^2}\right)dx}.
    \label{eq:conditional_density}
\end{equation}

\subsection{Conditional Normalization Constant}

The marginal density $f(\mathbf{x}_{2:d})$, which acts as the normalization constant of the conditional density in \eqref{eq:conditional_density}, is given by the integral
\begin{equation*}
f(\mathbf{x}_{2:d})=\int_{L}^{U} \, \phi\Bigl(\varepsilon \sqrt{(x - c_{1})^2 + \rho^2}\Bigr) \, dx.
\end{equation*}

To compute it, we apply a change of variables $u = x - c_{1}$, transforming the limits to $u \in [L- c_{1}, U- c_{1}]$. The integrand involves terms of the form $\sqrt{u^2 + \rho^2}$, which suggests a hyperbolic change of variables.
Let:
\[
    u = \rho \sinh(t) \implies du = \rho \cosh(t) dt, \quad \sqrt{u^2 + \rho^2} = \rho\cosh(t).
\]
Let us define the function 
\begin{equation}
\label{eq:tildex}
    \Tilde{X}(x):= \operatorname{arcsinh}\left(\frac{x-c_1}{\rho}\right)
\end{equation}
Then we can define the integration boundaries in the new variable $t$ as:
\[
    \Tilde{L} = \Tilde{X}(L), \quad \Tilde{U} = \Tilde{X}(U).
\]

The integrand involves the shape parameter $\varepsilon$, which is included directly in the primitive. Expressed in the hyperbolic variable $t$, we write
\begin{equation}
\label{eq:H0}
H_0(t)=\int_0^{\rho \sinh t}
\left(1-\varepsilon\sqrt{u^2+\rho^2}\right)^4\Bigl(4\varepsilon\sqrt{u^2+\rho^2}+1\Bigr)\,du,
\end{equation}
which can be computed explicitly and expressed in closed form in terms of hyperbolic functions.
\begin{align*}
H_0(t) &= \frac{1}{48}\Bigl[
48\rho\sinh t
+ \varepsilon^2\rho^3\bigl(-360\sinh t - 40\sinh 3t\bigr)\\
&+ \varepsilon^3\rho^4\bigl(360t + 240\sinh 2t + 30\sinh 4t\bigr)\\
&+ \varepsilon^4\rho^5\bigl(-450\sinh t - 75\sinh 3t - 9\sinh 5t\bigr)\\
&+ \varepsilon^5\rho^6\bigl(60t + 45\sinh 2t + 9\sinh 4t + \sinh 6t\bigr)
\Bigr].
\end{align*}

Applying the previous changes of variables,
\begin{equation}\label{eq:fx2d}
f(\mathbf{x}_{2:d})=\int_{L}^{U} \, \phi\Bigl(\varepsilon \sqrt{(x - c_{1})^2 + \rho^2}\Bigr) \, dx= H_0(\Tilde{U})-H_0(\Tilde{L}).
\end{equation}

Note that, unlike the univariate setting of Section \ref{sec:unipdf}, no case distinction on the position of the center $c_1$ is needed, since the effective bounds defined in \eqref{eq:LU} already account for the intersection of the component's support with $\Omega$. Moreover, this formulation is consistent with the univariate case: as $\rho \to 0$, the contribution of the conditioning variables vanishes, and therefore, one can verify that $H_0(\Tilde{U})-H_0(\Tilde{L})$ converges to the univariate normalization constant defined in \eqref{eq:normalization1d}, recovering directly the different cases relative to the position of the center.

\subsection{Conditional Moments}

As in Section \ref{sec:unipdf}, to obtain the actual normalized moments of the conditional density, we apply the binomial expansion and divide by the marginal density $f(\mathbf{x}_{2:d})$ (which acts as our conditional normalization constant). This yields the following expression for the moments:
\begin{align*}
    m_{k} &=  \frac{1}{f(\mathbf{x}_{2:d})} \int_{L}^{U} x^k \, \phi\Bigl(\varepsilon \sqrt{(x - c_{1})^2 + \rho^2}\Bigr) \, dx \\&=  \frac{1}{f(\mathbf{x}_{2:d})} \sum_{j=0}^{k} \binom{k}{j} c_{1}^{k-j} M_j,
\end{align*} 
where 
\begin{equation*}
    M_j= \int_{L-c_1}^{U-c_1}  u^j \phi(\varepsilon \sqrt{u^2+\rho^2})\,du.
\end{equation*}

Therefore, the computation of the conditional moments reduces to the computation of $M_j$. We will detail the computation for $j=1,2$, although it is computationally feasible to find primitives for higher-order terms and therefore to compute moments of much higher order.

\paragraph{Case $j=1$}
In this particular case, the factor $u = x-c_1$ in the integrand permits the change of variables $s^2 = \varepsilon^2(u^2+\rho^2)$, $2sds = 2\varepsilon^2
u\,du$, which eliminates the need for the hyperbolic substitution entirely. Thus, we can compute
\begin{equation*}
   M_1
    =\frac{1}{\varepsilon^2} \int_{r_L}^{r_U} s\phi(s) \, ds= \frac{1}{\varepsilon^2}\bigl[F_1(r_U) - F_1(r_L)\bigr],
    \label{eq:cond_m1}
\end{equation*}
where $F_1$ is the polynomial primitive defined in \eqref{eq:Fk} with $k=1$, and
\begin{equation*}
    r_L = \varepsilon\sqrt{(L-c_1)^2+\rho^2}, \qquad
    r_U = \varepsilon\sqrt{(U-c_1)^2+\rho^2}.
\end{equation*}

This expression is valid regardless of the position of $c_1$ relative to $[L,U]$, since the contributions of the two symmetric parts around $u=0$ cancel appropriately.

\paragraph{Case $j=2$}
For the second moment, that is $j=2$, we must compute the primitive
\begin{equation*}
H_2(t)=\int_0^{\rho \sinh t} u^2
\left(1-\varepsilon\sqrt{u^2+\rho^2}\right)^4\Bigl(4\varepsilon\sqrt{u^2+\rho^2}+1\Bigr)\,du,
\end{equation*}

which is given by the following hyperbolic expression:

\begin{align*}
   H_2(t) &= \frac{1}{5376}\Bigl[
\rho^3\bigl(-1344\sinh t + 448\sinh 3t\bigr)\\
&+ \varepsilon^2\rho^5\bigl(6720\sinh t - 1120\sinh 3t - 672\sinh 5t\bigr)\\
&+ \varepsilon^3\rho^6\bigl(-6720t - 1680\sinh 2t + 1680\sinh 4t + 560\sinh 6t\bigr)\\
&+ \varepsilon^4\rho^7\bigl(6300\sinh t - 420\sinh 3t - 756\sinh 5t - 180\sinh 7t\bigr)\\
&+ \varepsilon^5\rho^8\bigl(-840t - 336\sinh 2t + 168\sinh 4t + 112\sinh 6t + 21\sinh 8t\bigr) \Bigr].
\end{align*}

Using these primitives we can compute
\begin{equation*}
    M_2 = H_2(\Tilde{U}) - H_2(\Tilde{L}).
\end{equation*}

Once we compute these quantities, the computation of the moments is straightforward. As in the normalization case, we can recover the univariate moments when the conditioning variables vanish, that is, when $\rho \to 0$.

\subsection{Cumulative Distribution Function}
The conditional CDF is defined as 

\begin{equation*}
    F(x \mid \mathbf{x}_{2:d})=\int_{-\infty}^x f(t\mid \mathbf{x}_{2:d})\,dt,
\end{equation*}

so it requires integrating the density from the lower boundary of $\Omega=[a,b]$. With the primitives derived above, specifically $H_0(t)$ \eqref{eq:H0}, we can evaluate the CDF exactly without numerical integration. 

Given a point $x \in [a,b]$, the integration region is the intersection of $[L, U]$ with the interval $[a, x]$. Using $\Tilde{X}(x)$ as defined in \eqref{eq:tildex}, the normalized conditional CDF is then:
\begin{equation}
\label{eq:cond_cdf}
     F(x \mid \mathbf{x}_{2:d}) = 
    \begin{cases} 
        0 & \text{if } x < L , \\
        \frac{1}{f(\mathbf{x}_{2:d})}\left(H_0(\Tilde{X}(x)) - H_0(\Tilde{L})\right) & \text{if } x \in [L,U] , \\
        1 & \text{if } x >U.
    \end{cases}
\end{equation}

Note that the CDF is normalized by the marginal density $f(\mathbf{x}_{2:d})$, that is, the normalization constant that depends on the conditioning variables, ensuring the conditional CDF has value 1 at the upper boundary $U$.

\section{Mixtures of CS-RBFs and density estimation}\label{sec:mixtures}

 Finite mixture models (FMMs) \citep{Mclachlan2000} have established themselves as fundamental tools in statistical analysis, appearing in fields ranging from genetics and bioinformatics to engineering and economics. As stated by \cite{McLachan2019}, mixture models provide a semi-parametric framework for modeling unknown probability distributions. 

Formally, let $\mathbf{X} \in \mathbb{R}^d$ denote a $d$-dimensional random vector. A finite mixture density $f(\mathbf{X}; \boldsymbol{\theta})$ with $K$ components is defined as the weighted sum:
\begin{equation*}
    f(\mathbf{X}; \boldsymbol{\theta}) = \sum_{j=1}^{K} w_j \, f_j(\mathbf{X}; \boldsymbol{\theta}_j),
\end{equation*}
where $w_j$ are non-negative mixing weights satisfying $\sum_{j=1}^K w_j = 1$, and $f_j(\cdot)$ are the individual component densities of the mixture. 

Currently, the fitting of such models predominantly relies on the Expectation-Maximization (EM) algorithm to compute Maximum Likelihood (ML) estimates \citep{McLachan2019}. Historically, Gaussian Mixture Models (GMMs) \citep{Rasmussen1999} have been the default choice due to their simplicity and computational efficiency. 

However, real-world data frequently violate normality assumptions, exhibiting significant asymmetry, multi-modality, or heavy tails. To address these limitations, the literature has evolved towards non-normal mixture models. Significant developments include Skew-Normal or Skew-t Mixtures \citep{McLachan2019}, which offer more flexible, asymmetric mixture component shapes. 

Most of these models share one feature: they have an infinite support in $\mathbb{R}^d$. In this work, we address this limitation by employing mixtures of Compactly Supported Radial Basis Functions. While our individual kernels are isotropic, their localized nature allows the mixture to adapt to complex, asymmetric, and multi-modal data through the arrangement of components, while simultaneously ensuring the density is confined to a particular domain.

\subsection{Mixtures of CS-RBFs}
We define the CS-RBF mixture by combining the finite mixture framework with the Wendland CS-RBFs. We give a general formulation for conditional, univariate and multivariate density functions. In these cases, the probability density function is modeled as a weighted sum of $K$ spatially localized CS-RBFs. Our approach differs from the standard finite mixture frameworks in the normalization. We give a general normalization of the whole mixture instead of normalizing each component individually.

A univariate mixture density is formulated as:
\begin{equation*} 
    f(x) = \frac{1}{Z} \sum_{j=1}^K w_j \, \phi\bigl(\varepsilon_j |x-c_j|\bigr), 
\end{equation*}
with $c_j \in \mathbb{R}$, and mixing weights $w_j >0$. The normalization is given by $Z = \sum_{j=1}^K w_j N_j$, the weighted sum of the mass of each component, where $N_j$ is given by Equation~\eqref{eq:normalization1d} with the corresponding $\varepsilon_j$ and $c_j$. If we multiply and divide each component by its mass $N_j$ we get the form
\begin{equation*} 
    f(x) =  \sum_{j=1}^K \gamma_j \, f_j(x), \quad \gamma_j=\frac{w_j N_j}{Z},
\end{equation*}
where $f_j(x)=\phi(\varepsilon_j |x-c_j|)/N_j$ denotes the $j-$th individually normalized CS-RBF.

A useful property of this mixture formulation is the notion of responsibility weights. Given an observation $x_i$, the responsibility of the $j$-th component is defined as
\begin{equation*}
    r_{ij} = \frac{w_j\phi\bigl(\varepsilon_j |x_i - c_j|\bigr)}{\sum_{\ell=1}^K w_\ell \phi\bigl(\varepsilon_\ell |x_i - c_\ell|\bigr)}, 
 \end{equation*}
which is the probability that $x_i$ was generated by the $j$-th component. Unlike in standard Gaussian mixtures, the compact support of each CS-RBF means that $\phi(\varepsilon_j|x_i - c_j|) = 0$ whenever $x_i$ lies outside the support of component $j$. Therefore, typically many responsibilities $r_{ij}$ are exactly zero for any given $x_i$. This sparsity is a direct consequence of compact support.

Similarly, given $\mathbf{x} = (x_1, \dots, x_d) \in \mathbb{R}^d$, we formulate a conditional mixture as
\begin{equation*} 
    f(x_1|\mathbf{x}_{2:d}) =\frac{1}{Z(\mathbf{x}_{2:d})}  \sum_{j=1}^K w_j \,  \phi\bigl(\varepsilon_j \sqrt{(x_1-c_{1,j})^2+\rho_j^2}\bigr), 
\end{equation*}
with $\rho_j$ defined as in \eqref{eq:rho} for the $j$-th center $\mathbf{c}_j$, and $Z(\mathbf{x}_{2:d})=\sum_{j=1}^K w_j f_j(\mathbf{x}_{2:d})$ is the global normalization. Here $f_j(\mathbf{x}_{2:d})$ denotes the mass of the $j-$th component of the mixture, defined as in \eqref{eq:fx2d}. 

Note that we normalize the mixture at the end and not individually by component; otherwise we would be assigning each component a mass depending only on the general weights $w_j$ and not on the values where we condition.

Multiplying and dividing each term of the sum by its own mass $f_j(\mathbf{x}_{2:d})$ recovers the form 
\begin{equation*}
    f(x_1\mid\mathbf{x}_{2:d}) = \sum_{j=1}^K \gamma_j(\mathbf{x}_{2:d})\,f_j(x_1\mid\mathbf{x}_{2:d}),
\end{equation*}
where $f_j(x_1\mid\mathbf{x}_{2:d})$ denotes a normalized conditional density defined as in Section \ref{sec:condpdf} and 
\begin{equation*}
    \gamma_j(\mathbf{x}_{2:d}) = \frac{w_j\, f_j(\mathbf{x}_{2:d})}{Z(\mathbf{x}_{2:d})}, 
\end{equation*}
which is available in closed form and acts as a new weight depending on $\mathbf{x}_{2:d}$. Also note that, if the conditioning values fall outside of the support of a certain component, then $\gamma_j=0$ automatically, due to the compact support.
 
Because in both cases we have been able to write the mixture as a weighted linear combination of its components, all analytical properties derived in Sections \ref{sec:unipdf} and \ref{sec:condpdf} are preserved under the CS-RBF mixture structure. For example, given a $K$-component univariate mixture, its exact $k$-th moment is directly evaluated as
\begin{equation*}
    m_k = \sum_{j=1}^K \gamma_j m_{j,k},
\end{equation*}
where $m_{j,k}$ denotes the $k$-th moment of the $j$-th CS-RBF component, given in Equations \eqref{eq:moment} and \eqref{eq:centered_moments}; the same direct linear combination holds for the univariate CDF given by \eqref{eq:cdf_global} and \eqref{eq:cdf_outside}. In the conditional case, the analogous properties, including the conditional CDF \eqref{eq:cond_cdf}, are weighted sums under the new weights $\gamma_j(\mathbf{x}_{2:d})$.

In this same framework, we may formulate a joint density mixture as 
\begin{equation*}
    f(\mathbf{x})=\frac{1}{Z_d} \sum_{j=1}^K w_j \phi(\varepsilon_j \|\mathbf{x}-\mathbf{c}_j\|) ,
\end{equation*}
where the multivariate normalization constant is given by $Z_d=\sum_{j=1}^K w_j Z_{d,j}$, and
\begin{equation}
\label{eq:Nd}
    Z_{d,j}=\int_{\mathbb{R}^d} \phi\bigl(\varepsilon_j\|\mathbf{x}-\mathbf{c}_j\|\bigr)\,d\mathbf{x}
    = \frac{240\,\pi^{\frac{d}{2}}(d+1)\Gamma(d)}{(d+5)\,\Gamma(d/2)\Gamma(d+5)}\,\frac{1}{\varepsilon_j^d},
\end{equation}
are the normalization constants of the components $\phi(\varepsilon_j \| \mathbf{x} - \mathbf{c}_j \|)$ in $\mathbb{R}^d$.

This global normalization constant $Z_d$ is valid if we assume that the density is defined on the unbounded domain $\mathbb{R}^d$. However, if we consider any of the dimensions to be bounded, the normalization constant must instead be computed as the integral of
the density function restricted to that domain $\Omega \subset \mathbb{R}^d$. This would be the multivariate analogue of the truncation analyzed in Section~\ref{sec:unipdf}, which in general may not have a closed-form expression for $d>1$.

Consequently, in the joint multivariate case the treatment of bounded domains is lost, although the compact support of the components still confines any leakage to a neighborhood of the data. Exact conditional densities on bounded intervals can nevertheless be recovered from the fitted joint model through slicing and normalization. 

\subsection{Density estimation} \label{csrbfs:densityestimation}    

A natural application of a parametric family of densities is the estimation of a density function from data. The same overall procedure applies whether the target is a univariate density, or a conditional or multivariate density given multivariate data, differing only in the objective minimized.
              
We employ a constructive incremental algorithm that builds the mixture, allowing the model to choose its number of CS-RBF components. Initialization begins with a one-component mixture, where the weights and shape parameters are optimized via SGD. Our proposed algorithm iteratively increases the number of components, finding in each iteration a new set of centers by clustering and performing gradient descent, evaluating each candidate model on the training set; the held-out test set is used exclusively for the final reported metrics. To ensure computational efficiency, we implement an early stopping criterion that terminates the process if the chosen performance metric fails to improve for two consecutive iterations.  

The natural parameters of the mixture model are subject to strict constraints: the shape parameters and the mixing weights must be strictly positive ($\varepsilon_j > 0, w_j > 0$), and although it is not necessary for our mixture model, we also force the weights to sum to one ($\sum_{j=1}^K w_j = 1$). To optimize these parameters using standard unconstrained SGD, we introduce smooth, invertible transformations for these parameters \citep[Chap.~16]{Nocedal2006}. 

We then must find an optimal parameter configuration, which, as the centers are fixed in each iteration, reduces to estimating a vector $\boldsymbol{\theta}$ containing the transformed shape parameters, and mixture weights. We optimize these transformed shape and weight parameters by minimizing the negative log-likelihood with SGD \citep{Bottou2010,Bottou2018}. 

\textbf{Shape Parameters:} 
To ensure positivity, we map the shape parameters to the logarithmic domain by defining new unconstrained variables $\tilde{\varepsilon}_j = \log(\varepsilon_j)$, with the inverse mapping $\varepsilon_j = \exp(\tilde{\varepsilon}_j)$. 

\textbf{Weights:} 
To satisfy the weights constraint, the weights are written in terms of auxiliary parameters $v_j$, which are related to $w_j$'s via the softmax transformation:
\begin{equation*}
    w_j = \frac{e^{v_j}}{\sum_{k=1}^K e^{v_k}}.
\end{equation*}

\textbf{Centers:} In our framework, we fix the centers via a clustering algorithm so that they are not included during the gradient descent phase but are fixed at the beginning of each stage of the algorithm. In particular, we employ the k-means++ clustering algorithm \citep{Arthur2007}, which has an extremely low computational cost. For brevity, we refer to it simply as ``k-means'' throughout the remainder of the paper. Excluding the centers from the gradient descent phase has proven beneficial in practice, allowing a better optimization of the remaining parameters of the model.

We first split the sample into a training and a test set. We then employ as an objective function the Regularized Negative Log-Likelihood (NLL) with our new parameters over the training set. We define it as:
\begin{equation} \label{eq:objective}
    \mathcal{J}(\boldsymbol{\theta}) = - \frac{1}{M} \sum_{i=1}^M \log \left( \sum_{j=1}^K \frac{w_j}{Z} \phi\bigl(\varepsilon_j \|\mathbf{x}_i - \mathbf{c}_j\|\bigr) \right) + \lambda ( \|\mathbf{v}\|^2 + \|\tilde{\boldsymbol{\varepsilon}}\|^2),
\end{equation}
where $Z$ denotes the normalization constant, either the univariate one defined on \eqref{eq:normalization1d} or the multivariate version if $d>1$ defined on \eqref{eq:Nd}. Furthermore, $\lambda \ge 0$ is a hyperparameter for regularization. This term favors weight uniformity, that is, weights that do not differ much from each other, and penalizes excessively large or small shape parameters, thereby preventing overfitting. We choose an empirically determined regularization parameter which has proven to be effective in different density estimation tasks with different data. We define it as 
\[\lambda=M^{-0.8},\] 
which is dependent on the size of the dataset. This way, the regularization vanishes when the dataset grows.

During optimization, the gradients with respect to each of these new variables are computed explicitly via the chain rule. It is also worth noting that $Z_d$ depends on the dimension, and also on the weights and shape parameters and therefore must be included in the gradient computations.

\begin{algorithm} []
\caption{Density Estimation with CS-RBF mixture}
\label{alg:sgd_CS-RBF}
\begin{algorithmic}[1]
\State \textbf{Input:} Data $X$, Max Components $K$, Regularization $\lambda$, Epochs $T$, Batch $B$.
\State $\text{best\_metric} \leftarrow +\infty$
\State $i \leftarrow 0$ 

\For{k $=1$ to $K$}
\State \textbf{Initialization (Stage 1):}
\State \quad $\mathbf{c}_{1:k} \leftarrow \text{k-means}(X, k)$ \Comment{Fix centers}
\State \quad $v_{1:k} \leftarrow 1/k$ \Comment{Initialize equal weights}
\State \quad $\tilde{\varepsilon}_{1:k} \leftarrow \log(\mathcal{U}(0,1))$ \Comment{Random shape parameters}
\State \textbf{Optimization (Stage 2):}
\For{$t=1$ to $T$}
    \State $\eta_t \leftarrow \eta_{\min} + \frac{1}{2}(\eta_{\max} - \eta_{\min})\bigl(1 + \cos(\pi \tfrac{t-1}{T})\bigr)$
    \State Sample a random mini-batch $X_b \subset X$ of size $B$
    \State Compute gradient $\nabla_{\tilde{\boldsymbol{\theta}}} \mathcal{J}$ on $X_b$
    \State Update $\tilde{\boldsymbol{\theta}} \leftarrow \tilde{\boldsymbol{\theta}} - \eta_t \nabla_{\tilde{\boldsymbol{\theta}}} \mathcal{J}$
\EndFor
\State $\text{score}_k \leftarrow \text{Metric}(X, \boldsymbol{\theta}_k)$ \Comment{Evaluate NLL or BIC}    
    \If{$\text{score}_k < \text{best\_metric}$}
        \State $\text{best\_metric} \leftarrow \text{score}_k$
        \State $\boldsymbol{\theta}_{\text{best}} \leftarrow \boldsymbol{\theta}_k$
        \State $\text{i} \leftarrow 0$
    \Else
        \State $\text{i} \leftarrow \text{i} + 1$
    \EndIf
    
    \If{$\text{i} \ge 2$}
        \State \textbf{break} \Comment{Early stop}
    \EndIf

\EndFor
\State \textbf{Output:} Transform $\tilde{\boldsymbol{\theta}}$ back to the constrained space and include the centers, weights, and shape parameters $\{\mathbf{c}, w, \varepsilon\}$.
\end{algorithmic}
\end{algorithm}

Furthermore, our SGD framework requires that we choose a maximum learning rate, which we generally fix as $\eta_{max}=1$, and establish a minimum of $\eta_{min}=10^{-5}$. At each epoch of the SGD, the learning rate decays following a cosine annealing schedule, as proposed by \citet{Loschilov2017}; this schedule helps avoid poor local minima during optimization. We also need to choose a number of epochs for the algorithm, which we generally set as $T=250$. For each iteration of the algorithm, we initialize, for simplicity's sake, with equal weights and shape parameters drawn uniformly from the interval $[0,1]$. The detailed algorithm is given in Algorithm~\ref{alg:sgd_CS-RBF}.

In the conditional case, we follow the same procedure described in the univariate case to fit the conditional mixture of Section~\ref{sec:mixtures} directly against pairs $(x_1, \mathbf{x}_{2:d})$. The procedure shares the reparameterization, the k-means center initialization, and the incremental scheme of Algorithm~\ref{alg:sgd_CS-RBF}. Only the loss function changes from $\mathcal{J}$ to the Conditional Negative Log-Likelihood (CNLL):
\begin{equation} \label{eq:objective_cond}
\begin{split}
    \mathcal{J}_{\text{cond}}(\boldsymbol{\theta}) = {}& -\frac{1}{M}\sum_{i=1}^M \log\left(\frac{1}{Z(\mathbf{x}_{2:d,i})}\sum_{j=1}^K w_j\,\phi\bigl(\varepsilon_j\sqrt{(x_{1,i}-c_{1,j})^2+\rho_{j,i}^2}\bigr)\right) \\
    & + \lambda\bigl(\|\mathbf{v}\|^2+\|\tilde{\boldsymbol\varepsilon}\|^2\bigr),
\end{split}
\end{equation}
with $\rho_{j,i}$ and $Z(\mathbf{x}_{2:d,i})$ as defined in Section~\ref{sec:mixtures}, and $\lambda$ chosen exactly as in \eqref{eq:objective}. Algorithm~\ref{alg:sgd_CS-RBF} applies unchanged with $\mathcal{J}_{\text{cond}}$.

\section{Experimental Results}\label{sec:experiments}
To demonstrate the ability of CS-RBFs to model density functions, we evaluate the empirical performance of the CS-RBF mixture framework for density estimation. We test it on both synthetic and real-world datasets, estimating univariate, multivariate, and conditional density functions.  

To measure the model's fit, we use the Bayesian Information Criterion (BIC) \citep{Schwartz1978}, defined as follows:
\begin{equation*}
    \text{BIC} = P \ln(M) - 2 \ln(\hat{L}),
\end{equation*} 
where $P$ is the number of estimated parameters, $M$ is the number of data points, and $\hat{L}$ is the maximized likelihood of the model on the data. A lower BIC indicates a better model in terms of likelihood of the data and model simplicity. In our experiments, the BIC is computed on the held-out test set, using the number of
free parameters of each model and the model's true, unpenalized log-likelihood, not the regularized training objective $\mathcal{J}$ of \eqref{eq:objective}.

We compare our model against a standard baseline: a (non-truncated) Gaussian Mixture Model (GMM) with full covariance matrices and regularization. In every experiment, the number of GMM components is selected in each run as the value minimizing the BIC over a range of candidate values, with the same early stopping rule as our incremental algorithm; reported training times include this model selection stage for both methods.

For a CS-RBF mixture with $K$ components in dimension $d$, the number of free parameters is $P = K(d+2)-1$: $K-1$ mixing weights, since they are subject to the sum-to-one constraint, $K$ shape parameters, and $Kd$ center coordinates. The centers are included in the count even though they are fixed by k-means, since they are estimated from the data.

For a GMM with full covariance matrices, $P = K\bigl(1+d+\frac{d(d+1)}{2}\bigr)-1$. In the univariate case, both counts coincide ($3K-1$); in higher dimensions, each Gaussian component spends $d(d+1)/2$ parameters on its covariance matrix, whereas each isotropic CS-RBF component uses a single shape parameter.

Furthermore, we define a ``boundary leak'' metric: the total probability mass assigned outside a previously known valid interval $[a,b]$, computed by numerically integrating the fitted density outside the variable domain. This leak is exactly $0$ for our model, which analytically self-normalizes within bounded domains. For the synthetic experiments, where the true density is known, we additionally report the mean squared error (MSE) between the estimated and the true density, evaluated on a uniform grid of 1000 points on an extended interval containing the variable domain, so as to also account for the error due to the boundary leak.

All reported results are averaged over multiple independent runs, to account for the stochastic nature of both the sampling and the fitting algorithms. For the number of components $K$ we also report the standard deviation across runs, since it reflects the stability of the model selection.

\subsection{Univariate Estimation}
First, we evaluate the estimation of univariate probability density functions. For synthetic data, we compare both a standard regularized GMM and our CS-RBF mixture against a known ground-truth density. For real-world datasets, the two models are compared against each other using standard information criteria. We also compute the first moments of the models and the known distributions to measure the error introduced by the models.

\textbf{Beta mixture}

The first experiment we perform is to fit a density function using artificially generated samples from a known distribution. We generate $M=2000$ random samples from a mixture of two Beta distributions, $0.6\,\mathcal{B}(2,12) + 0.4\,\mathcal{B}(12,2)$, explicitly chosen to create a bounded ($x \in [0,1]$), multimodal, and skewed density. We split the data into a training set containing $80\%$ of the samples and a test set used to evaluate the performance of the model. We measure the BIC and the boundary leak metric, in this case the variable domain is $[0,1]$.

\begin{table}[htbp]
  \centering
  \caption{Performance metrics averaged over 50 runs.}
  \label{tab:exp11}
  \begin{tabular}{lcc}
    \hline
    Metric & CS-RBF & GMM \\
    \hline
    Log-Likelihood & 169.48 & 167.32 \\
    BIC        & -234.59  & -239.25  \\
    MSE         & 0.031130 & 0.027350  \\
    Boundary Leak & 0 & 0.008257 \\
    Training Time (s)         & 0.1392 & 0.2606  \\
    Num.\ Components ($K$)   & 6.1 $\pm$ 0.7 & 5.6 $\pm$ 0.6 \\
    \hline
  \end{tabular}
\end{table}

Table \ref{tab:exp11} summarizes the performance of the two models over 50 runs with different seeds. Both models approximate the true density accurately: the CS-RBF mixture attains a slightly higher test log-likelihood, while the GMM, selecting slightly fewer components, obtains a marginally better MSE and BIC. In this experiment the GMM leaks only about $0.8\%$ of its probability mass outside the domain. As the following experiments show, the advantage of the CS-RBF model grows precisely when the true density is positive at a boundary of the domain. The training times of both models, which include model selection, are of the same order of magnitude.

\begin{table}[htbp]
  \centering
  \caption{Relative moment errors with respect to the true distribution averaged over 50 runs.}
  \label{tab:exp12}
  \begin{tabular}{ccc}
    \hline
    Moment & CS-RBF (\%) & GMM (\%) \\
    \hline
    $m_1$ & 1.71  & 0.85  \\
    $m_2$ & 2.41  & 1.17  \\
    $m_3$ & 2.69  & 1.32  \\
    $m_4$ & 2.88  & 1.46  \\
    \hline
  \end{tabular}
\end{table}

Table \ref{tab:exp12} reports the relative error of each model's moments with respect to the true moments of the underlying distribution. The GMM consistently yields smaller moment errors than CS-RBF. Our algorithm does not enforce the mixture to reproduce the empirical moments, since it only tries to minimize negative log-likelihood. Even so, the CS-RBF model still reproduces the theoretical moments with a low percentage of error, below $3\%$.

\begin{figure}[]
    \centering
\includegraphics[width=0.7\linewidth]{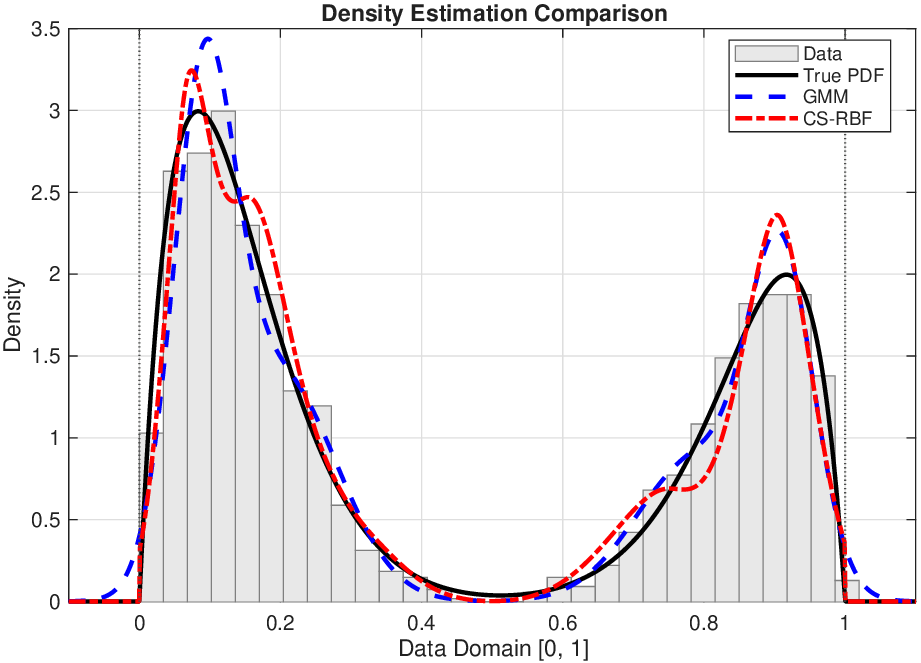}
    \caption{Comparison of the estimated density functions for the synthetic bimodal Beta distribution. The GMM leaks probability below $0$ and above $1$, whereas the CS-RBF strictly respects the boundaries.}
    \label{fig:beta_density}
\end{figure}

Figure~\ref{fig:beta_density} illustrates the estimated densities of the model: the GMM's density curve remains visibly non-zero outside $[0,1]$, while the CS-RBF's is truncated exactly at the boundary by construction.

\textbf{Pareto distribution}

In this second experiment, we generate $M=1000$ random samples from a Generalized Pareto distribution (GPD) with shape $\kappa=-0.25$, scale $\sigma=3$, and location $\mu=1$, whose support is the bounded interval $[\mu,\,\mu-\sigma/\kappa]=[1,13]$. 

We split the data into training and test sets. Table \ref{tab:exp21} summarizes the evaluation metrics. Unlike the previous experiment, the true density is strictly positive at the lower endpoint $x=\mu$. The GMM leaks about $1.2\%$ of the probability mass outside the domain, and the model requires roughly twice as many components ($K \approx 4$) as our mixture ($K \approx 2$) to represent the sharp truncation. As a result, the CS-RBF model attains a clearly better BIC at essentially identical log-likelihood. The GMM tracks the interior shape of the distribution slightly better, as reflected in its lower MSE and CDF error, but at the cost of the leaked mass and the additional components. By contrast, our Wendland $\mathscr{C}^2$ mixture leverages its compact support and closed-form normalization to eliminate any boundary bias.

\begin{table}[htbp]
  \centering
  \caption{Performance metrics for the GPD experiment averaged over 50 runs.}
  \label{tab:exp21}
  \begin{tabular}{lcc}
    \hline
    Metric & CS-RBF & GMM \\
    \hline
    Log-Likelihood   & -374.77 & -374.76 \\
    BIC         & 777.94  & 809.72  \\
    MSE        & 0.000847  & 0.000616  \\
    CDF MSE      & 0.001010  & 0.000118 \\
    Boundary Leak   & 0 & 0.011888 \\
    Training Time (s)           & 0.0747 & 0.1163  \\
    Num.\ Components ($K$)      & 2.1 $\pm$ 0.3 & 4.1 $\pm$ 0.6 \\
    \hline
  \end{tabular}
\end{table}

\begin{figure}[]
    \centering
    \includegraphics[width=1\linewidth]{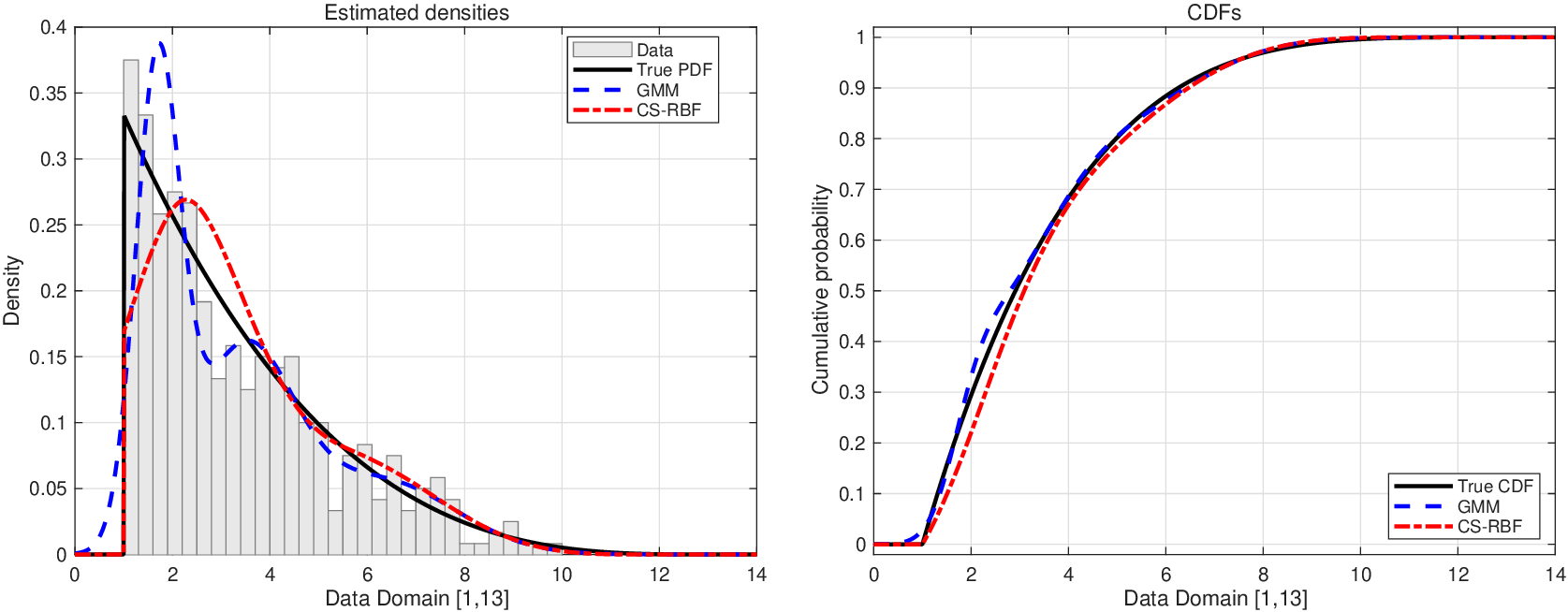}
    \caption{Comparison of the estimated density functions for the synthetic Generalized Pareto distribution (left). Cumulative distribution functions of the two models compared to the true CDF of the distribution (right).}
    \label{fig:gpd_density}
\end{figure}

Figure~\ref{fig:gpd_density} shows in the left panel the estimated densities together with the true density of the distribution: our model is truncated exactly at the lower bound of the variable and vanishes naturally at the upper bound, whereas the GMM leaks probability below $x=1$. The right panel shows the estimated CDFs together with the theoretical one. The corresponding mean squared CDF errors are reported in Table~\ref{tab:exp21}.

\textbf{Real data}

To validate our model on real physical measurements, we perform a third experiment using real data from the widely known \textit{Abalone} dataset \citep{abalone}, selecting four continuous variables with distinct distributional characteristics.

\begin{figure}[]
    \centering
\includegraphics[width=0.8\linewidth]{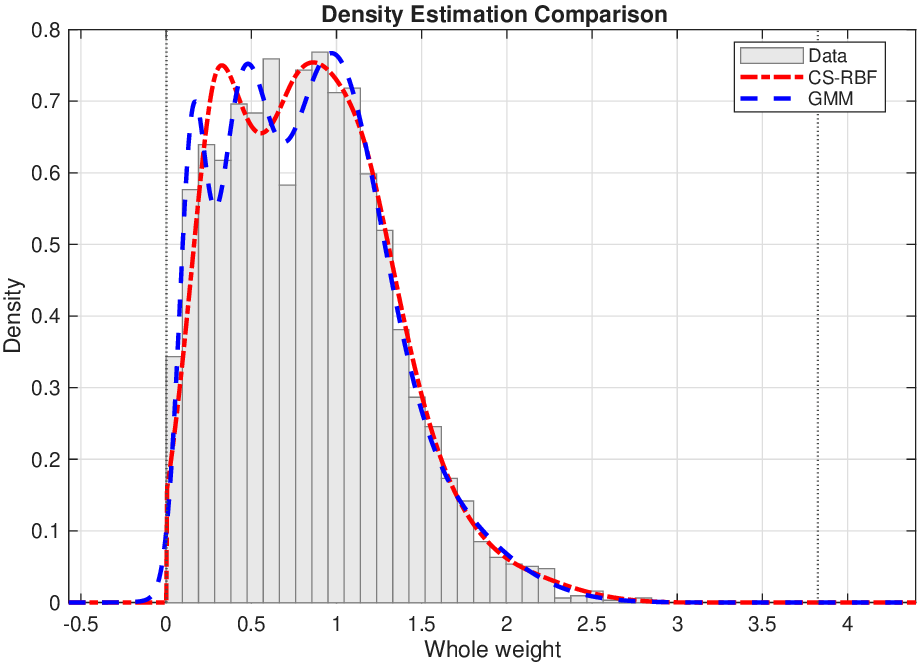}
    \caption{Probability density functions estimated by the CS-RBF model and a standard GMM for the \textit{Whole weight} variable.}
    \label{fig:densityreal}
\end{figure}

\begin{table}[htbp]
  \centering
  \caption{Performance metrics across Abalone variables
  averaged over 10 runs.}
  \label{tab:abalone_metrics}
  \begin{tabular}{lcc}
    \hline
    Metric & CS-RBF & GMM \\
    \hline
    \multicolumn{3}{l}{\textit{Length}} \\
    \hline
    BIC       & -1193.84  & -1212.18  \\
    Boundary Leak & 0 & 0.000302 \\
    Num.\ Components ($K$)    & 2.4 $\pm$ 0.5 & 2.0 $\pm$ 0.0 \\
    \hline
    \multicolumn{3}{l}{\textit{Height}} \\
    \hline
    BIC    & -1513.83 & -1531.92  \\
    Boundary Leak & 0 & 0.000508 \\
    Num.\ Components ($K$)    & 2.1 $\pm$ 0.3 & 2.0 $\pm$ 0.0 \\
    \hline
    \multicolumn{3}{l}{\textit{Whole weight}} \\
    \hline
    BIC      & 1143.45 & 1146.30  \\
    Boundary Leak & 0 & 0.004187  \\
    Num.\ Components ($K$)    & 4.0 $\pm$ 0.5 & 4.0 $\pm$ 0.0 \\
    \hline
    \multicolumn{3}{l}{\textit{Shell weight}} \\
    \hline
    BIC         & -972.52  & -948.22  \\
    Boundary Leak & 0 & 0.005331  \\
    Num.\ Components ($K$)    & 3.0 $\pm$ 0.0 & 4.2 $\pm$ 0.4 \\
    \hline
  \end{tabular}
\end{table}

Figure~\ref{fig:densityreal} illustrates the \textit{Whole weight} variable, a physically bounded quantity (positive) with a mode close to the boundary. The results in Table~\ref{tab:abalone_metrics} reveal a consistent trend across all four variables: as the GMM's boundary leak grows, the model is outperformed by our truncated CS-RBF model. Therefore, the CS-RBF model advantage is not fixed, but grows in step with how much probability mass a model like the GMM is forced to place in physically impossible regions. Both models have a comparable number of free parameters across the four variables, so the differences in BIC are attributable to fit and boundary behavior rather than to model size.

\subsection{Multivariate estimation}

As discussed previously, the CS-RBF mixture is naturally defined in multivariate spaces. However, some aspects, such as a closed-form normalization over arbitrary bounded domains, are only available in the univariate and conditional settings. To demonstrate the model's ability to estimate joint densities, we estimate the joint density of two variables from the \textit{Bodyfat} dataset \citep{Johnson1996}.

Since the normalization over the intersection of the components' supports with a bounded domain has no closed-form expression for $d>1$, the joint CS-RBF mixture is normalized with the unbounded constant of Equation~\eqref{eq:Nd}, and its practical validity is tested against the GMM baseline.
 
\begin{figure}[]
    \centering
    \includegraphics[width=0.95\linewidth]{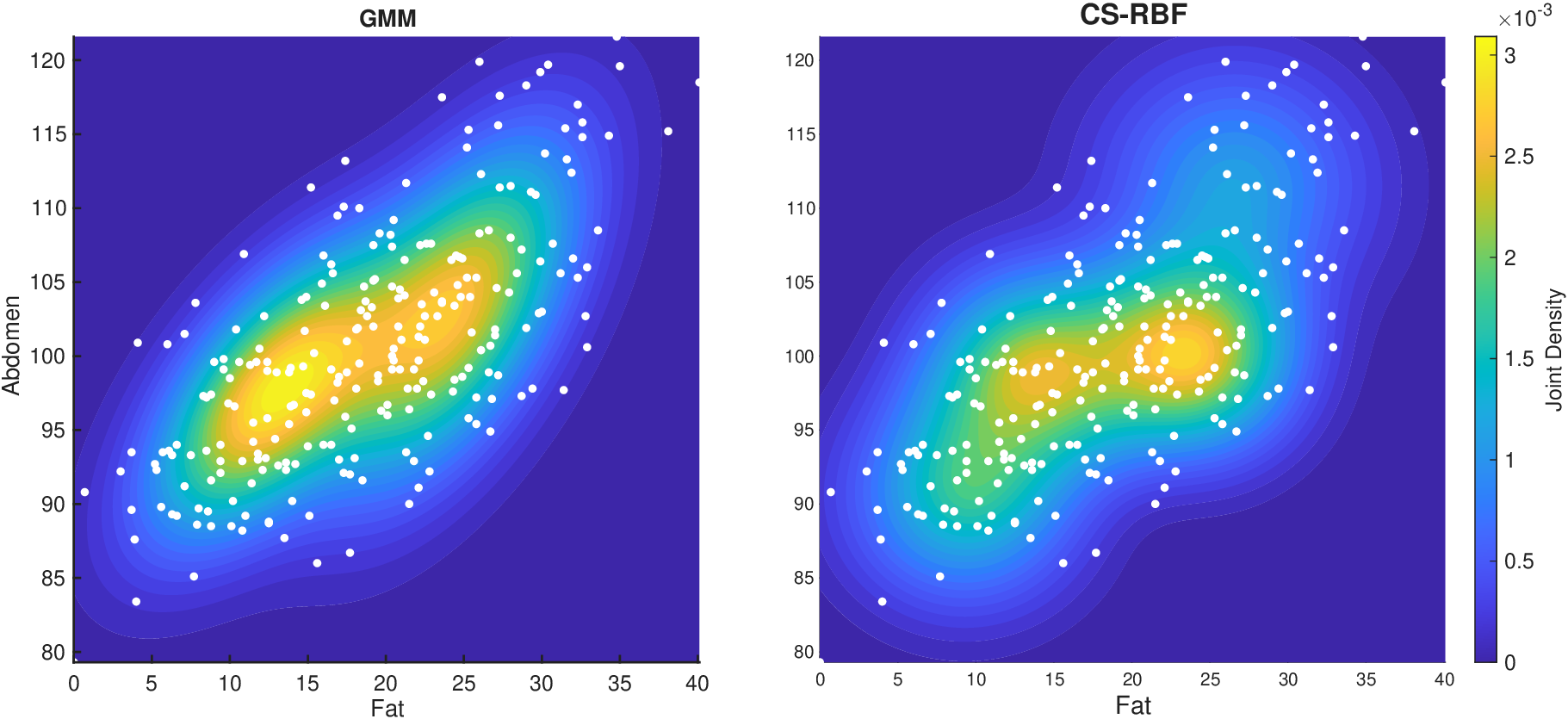}
    \caption{Joint density estimation of two variables from the
    \textit{Bodyfat} dataset with a GMM (left) and the CS-RBF model (right).}
    \label{fig:bivariate}
\end{figure}
 
As shown in Figure~\ref{fig:bivariate}, the proposed model estimates bivariate densities without requiring an excessive number of components. Table~\ref{tab:bivariate_metrics} reports the results for two variable pairs of the dataset, averaged over 10 runs. In both pairs the GMM attains a slightly higher log-likelihood and a better BIC: its full-covariance components capture the correlation between the variables directly, whereas the isotropic CS-RBF kernels must approximate that correlation through a larger number of components ($K\approx4$--$6$ versus $K=2$). The CS-RBF mixture nevertheless remains close in fit, showing that the model is a practical option for joint estimation. Moreover, while the joint estimation assumes an unbounded domain, the fitted bivariate density can subsequently be transformed into a strict conditional density on a bounded interval by slicing at a specific value of one variable and renormalizing analytically with the marginal derived in Section~\ref{sec:condpdf}.

\begin{table}[htbp]
  \centering
  \caption{Bivariate joint density estimation on two variable pairs of
  the \textit{Bodyfat} dataset averaged over 10 runs.}
  \label{tab:bivariate_metrics}
  \begin{tabular}{lcc}
    \hline
    Metric & CS-RBF & GMM \\
    \hline
    \multicolumn{3}{l}{\textit{Fat / Chest}} \\
    \hline
    Log-Likelihood  & -6.7872 & -6.7132 \\
    BIC & 726.63 & 700.70 \\
    Num.\ Components ($K$)  & 4.2 $\pm$ 0.8 & 2.0 $\pm$ 0.0 \\
    \hline
    \multicolumn{3}{l}{\textit{Weight / Abdomen}} \\
    \hline
    Log-Likelihood  & -7.9794 & -7.7095 \\
    BIC   & 880.09  & 813.98 \\
    Num.\ Components ($K$)  & 5.5 $\pm$ 0.5 & 2.0 $\pm$ 0.0 \\
    \hline
  \end{tabular}
\end{table}

\subsection{Conditional Density Estimation}

We now turn from the univariate case to conditional density estimation, using the direct training procedure of Section~\ref{csrbfs:densityestimation} and the CNLL objective in \eqref{eq:objective_cond}. To evaluate this capability, we use four configurations from three real datasets: \textit{Glass} \citep{glass}, \textit{Chickenpox} \citep{chickenpox} with two different sets of variables, and \textit{Bodyfat} \citep{Johnson1996}.

We compare a joint multivariate GMM model, from which conditionals are obtained, with our CS-RBF mixture. In each configuration, we select one target variable $x_1$ and three conditioning variables $\mathbf{x}_{2:4}$, forming a 4-dimensional joint space, and measure the Conditional Log-Likelihood (CLL) on test data, together with the BIC, the boundary leak per test point, the training time, and the number of components selected by each model. For the GMM, the number of components is selected in each run as the value minimizing the conditional BIC on the training set. The results can be seen in Table~\ref{tab:conditional_metrics}.

In each configuration, the target variable is, respectively: calcium content
(\textit{Glass}); weekly chickenpox cases in Borsod (\textit{Chickenpox}~(A));
abdomen circumference (\textit{Bodyfat}); and weekly cases in Budapest
(\textit{Chickenpox}~(B)). The three conditioning variables are chosen among
the remaining continuous variables of each dataset.

\begin{table}[htbp]
  \centering
  \caption{Conditional density estimation across four target/conditioning configurations averaged over 10 runs.}
  \label{tab:conditional_metrics}
  \begin{tabular}{lcc}
    \hline
    Metric & Cond.\ CS-RBF & Joint GMM \\
    \hline
    \multicolumn{3}{l}{\textit{Glass} } \\
    \hline
    CLL   & -3.96  & 8.52 \\
    BIC  & 88.54 & 161.93 \\
    Boundary leak & 0 & 0.0231 \\
    Training time (s)          & 0.2810  & 0.0268  \\
    Num.\ Components ($K$)     & 4.7 $\pm$ 1.3 & 4.1 $\pm$ 1.5 \\
    \hline
    \multicolumn{3}{l}{\textit{Chickenpox (A)} } \\
    \hline
    CLL       & -251.43 & -241.68\\
    BIC     & 581.74 & 597.73  \\
    Boundary leak  & 0 & 0.0152 \\
    Training time (s)          & 0.3606 & 0.0275  \\
    Num.\ Components ($K$)     & 3.5 $\pm$ 1.3 & 2.0 $\pm$ 0.0 \\
    \hline
    \multicolumn{3}{l}{\textit{Bodyfat}} \\
    \hline
    CLL  & -69.71  & -67.19  \\
    BIC     & 316.61  & 227.13  \\
    Boundary leak  & 0 & 0.0002 \\
    Training time (s)          & 0.6789 & 0.0305\\
    Num.\ Components ($K$)     & 9.4 $\pm$ 0.8 & 2.0 $\pm$ 0.0 \\
    \hline
    \multicolumn{3}{l}{\textit{Chickenpox (B)} } \\
    \hline
    CLL & -264.78  & -263.59 \\
    BIC   & 648.71  & 641.59  \\
    Boundary leak & 0 & 0.0078 \\
    Training time (s)          & 0.4448  & 0.0387 \\
    Num.\ Components ($K$)     & 5.2 $\pm$ 0.8 & 2.0 $\pm$ 0.0 \\
    \hline
  \end{tabular}
\end{table}

In general, the joint GMM attains slightly higher conditional log-likelihoods, as its
full-covariance components directly capture the linear dependence between the target and the conditioning variables. This gap is largest for \textit{Glass}, as can be seen in Table~\ref{tab:conditional_metrics}.

In terms of BIC, in the configurations with the largest boundary leak (\textit{Glass} and \textit{Chickenpox}~(A)), the CS-RBF mixture attains a clearly lower BIC with a comparable number of components, whereas in configurations where boundary effects are negligible (\textit{Bodyfat}), the GMM is more parsimonious and obtains the better BIC. The CS-RBF mixture assigns exactly zero mass outside the domain in every case. In \textit{Bodyfat}, the isotropic CS-RBF kernels require considerably more components ($K \approx 9$) to capture correlations than a full-covariance Gaussian model captures directly. This also explains its larger BIC in that specific setup.

Regarding training time, the GMM is roughly one order of magnitude faster than our CS-RBF model. This is partially due to the datasets' sample sizes, around $M=500$. Although the training times of our model remain far from prohibitive, we believe that the fitting algorithm of the conditional CS-RBF mixture can be further optimized.

\begin{figure}[]
    \centering
    \includegraphics[width=0.95\linewidth]{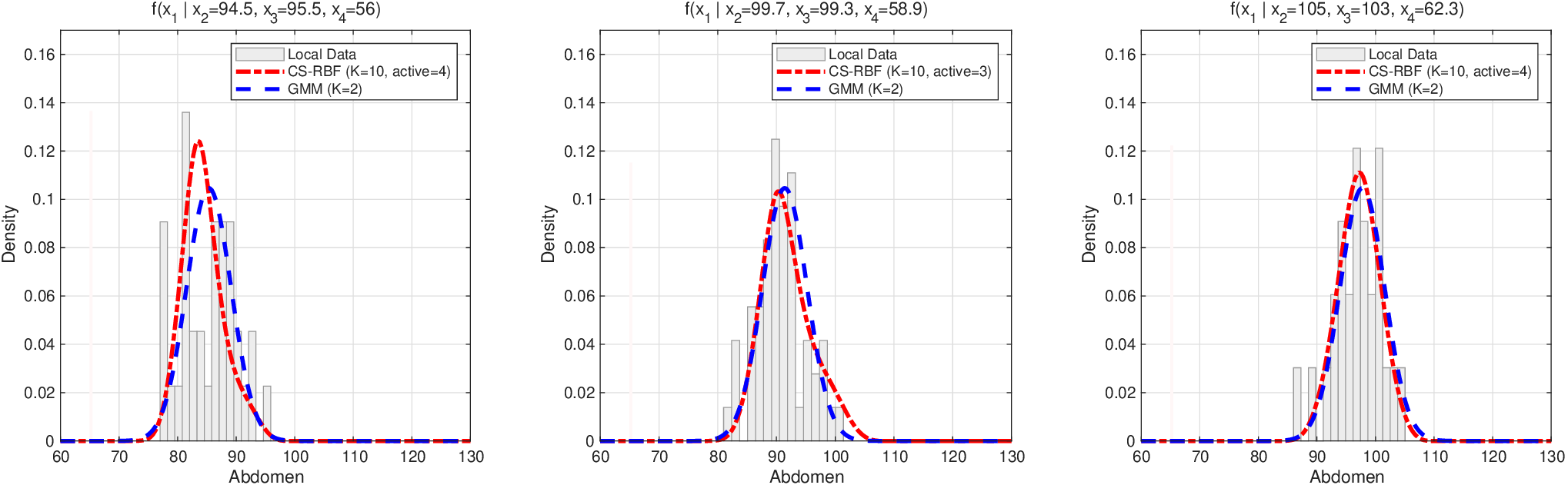}
    \caption{Conditional density slices for the \textit{Bodyfat} configuration, at three different configurations of the conditioning variables. Grey bars show the empirical histogram of the data points close to the conditional slice.}
    \label{fig:4}
\end{figure}

Figure \ref{fig:4} shows the conditional densities estimated by both models on the \textit{Bodyfat} dataset, sliced at three different configurations of the conditioning variables, corresponding to their quartiles. For each configuration, we also indicate which CS-RBFs of the mixture are active, that is, those that are nonzero at that conditioning point.

\section{Conclusion}

In this paper, we formalize Compactly Supported Radial Basis Functions (CS-RBFs) as a novel parametric family of probability density functions. For simplicity, the focus is on Wendland $\mathscr{C}^2$ kernels, but the approach can be extended to other families. The major contribution of this work is providing analytical closed-form expressions for the main statistical properties of CS-RBF densities. Both univariate and conditional densities are explored, and explicit expressions for normalization constants, moments, and cumulative distribution functions are obtained in each case. In each case, the results comprise the cases when the CS-RBF support is truncated, due to restrictions in the continuous variable domain, and the untruncated case. Having exact, closed-form analytical expressions is an advantage over other approaches that often rely on numerical approximation techniques.

On this foundation, we consider combinations of CS-RBFs to build mixture models, for which most properties are derived from those of individual CS-RBFs. CS-RBF mixture models, as in the case of Gaussian functions and GMMs, offer greater flexibility while keeping the model complexity controlled. As a second significant contribution, we develop a learning procedure for CS-RBF mixture models for density estimation, and we evaluate it experimentally. 

The results show that CS-RBF mixtures exhibit good performance in modeling synthetic and real-world data. This is particularly evident for continuous random variables with a bounded domain. These variable constraints can be naturally incorporated into the CS-RBF model. Due to their properties, CS-RBF mixtures yield competitive results in comparison to GMMs in terms of likelihood and BIC. The use of CS-RBF density mixtures has some extra benefits: all probability mass is assigned within the variable domain (no boundary bias), and the compact support of its components improves the efficiency of its evaluation and prevents numerical problems. These results highlight the effectiveness of CS-RBFs as an efficient tool to model continuous variables in density estimation, while enabling exact computation of their main probabilistic properties.

The approach presented in this work may be extended in various ways. On the one hand, the mixture learning might be improved by allowing CS-RBF centers to adapt dynamically or by exploring alternative formulations to the optimization process. While this paper established the foundational theory of CS-RBFs as density functions using Wendland $\mathscr{C}^2$ functions, the framework can be naturally extended to other CS-RBFs of higher smoothness. Furthermore, the current isotropic framework provides a robust baseline for transitioning to anisotropic CS-RBFs. By utilizing the Mahalanobis norm ($A$-norm), the model could more effectively capture directional dependencies in high-dimensional or highly correlated feature spaces, although the analytic derivation of the conditional densities after ``slicing'' becomes much more involved.

The presented closed-form expressions for univariate and conditional densities, including those with truncated support, allow the integration of CS-RBFs into more complex probabilistic models. Specifically, in the context of hierarchical models and Bayesian networks, where each node is naturally modeled by either a univariate or a conditional density function. In such models, the compact support of CS-RBFs and the possibility of exact computation of moments or the CDF could be further exploited for different purposes.

\section*{Acknowledgments}

All authors are supported by the project PID2025-170285NB-I00 funded by MICIU/AEI/10.13039/501100011033 and by ``ERDF A way of making Europe''.

The first author (SDE) is supported by the University of Almería’s program for research and knowledge transfer (PPIT-UAL).

The second author (AMF) was partially supported by Junta de Andaluc\'{\i}a (Instituto Interuniversitario Carlos I de F\'{\i}sica Te\'orica y Computacional). 

The third author (DRL) is also supported by the projects PID2022-139293NB-C31 and PID2025-172960NB-C33 funded by MICIU/AEI/10.13039/501100011033 and by ``ERDF A way of making Europe''.

The authors thank the support from the research group FQM-229 (Junta de Andalucía), and from the Center for Development and Transfer of Mathematical Research to Industry CDTIME (University of Almería).

\bibliographystyle{unsrtnat}
\bibliography{refs}  

@book{Silverman1986,
  title     = {Density Estimation for Statistics and Data Analysis},
  author    = {Silverman, Bernard W.},
  year      = {2018},
  publisher = {Routledge},
  address= {London}
}

@book{Wendland2005, address= {Cambridge}, series={Cambridge Monographs on Applied and Computational Mathematics}, title={Scattered Data Approximation}, publisher={Cambridge University Press}, author={Wendland, Holger}, year={2004}, collection={Cambridge Monographs on Applied and Computational Mathematics}}

@book{Scott2015,
author = {Scott, David},
year = {2015},
month = {03},
title = {Multivariate density estimation: Theory, practice, and visualization: Second edition},
isbn = {9780471697558},
publisher={John Wiley \& Sons},
doi = {10.1002/9781118575574},
address={New York}
}

@article{Jawitz2004,
  title   = {Moments of truncated continuous univariate distributions},
  author  = {Jawitz, James W.},
  journal = {Advances in Water Resources},
  volume  = {27},
  pages   = {269--281},
  year    = {2004},
  doi = {10.1016/j.advwatres.2003.12.002}
}

@inproceedings{Rasmussen1999,
author = {Rasmussen, Carl},
year = {2000},
month = {04},
pages = {554-560},
title = {The Infinite Gaussian Mixture Model},
volume = {12},
booktitle = {Adv Neural Inf Process Syst}
}

@article{Cattaneo2020,
  author    = {Cattaneo, Matias D. and Jansson, Michael and Ma, Xinwei},
  title     = {Simple Local Polynomial Density Estimators},
  journal   = {Journal of the American Statistical Association},
  volume    = {115},
  number    = {531},
  pages     = {1449--1455},
  year      = {2020}
}

@article{Cattaneo2022,
  author    = {Cattaneo, Matias D. and Jansson, Michael and Ma, Xinwei},
  title     = {{lpdensity}: {Local} Polynomial Density Estimation and Inference},
  journal   = {Journal of Statistical Software},
  volume    = {101},
  number    = {2},
  pages     = {1--25},
  year      = {2022}
}

@article{Mclachan2019,
author = {Mclachlan, G. and Lee, Sharon and Rathnayake, Suren},
year = {2019},
month = {03},
pages = {355-378},
title = {Finite Mixture Models},
volume = {6},
journal = {Annual Review of Statistics and Its Application},
doi = {10.1146/annurev-statistics-031017-100325}
}

@book{Fasshauer2007,
  title={Meshfree approximation methods with {MATLAB}},
  author={Fasshauer, Gregory E},
  year={2007},
  publisher={World Scientific Publishing Company},
  address = {New York}
}

@article{Baes2024,
  author  = {Baes, Maarten},
  title   = {Self-consistent dynamical models with a finite extent -- {IV}. {Wendland} models based on compactly supported radial basis functions},
  journal = {Monthly Notices of the Royal Astronomical Society},
  volume  = {531},
  number  = {4},
  pages   = {5097--5108},
  year    = {2024},
  doi     = {10.1093/mnras/stae1521}
}

@article{Dehnen2012,
    author = {Dehnen, Walter and Aly, Hossam},
    title = {Improving convergence in smoothed particle hydrodynamics simulations without pairing instability},
    journal = {Monthly Notices of the Royal Astronomical Society},
    volume = {425},
    number = {2},
    pages = {1068-1082},
    year = {2012},
    month = {09},
    issn = {0035-8711},
    doi = {10.1111/j.1365-2966.2012.21439.x},
}

@book{Bishop2006,
  title={Pattern recognition and machine learning},
  author={Bishop, Christopher M},
  year={2006},
  publisher={Springer},
  address = {New York}
}

@article{Colbrook2020,
  author    = {Colbrook, Matthew J. and Botev, Zdravko I. and Kuritz, Karsten and MacNamara, Shev},
  title     = {Kernel Density Estimation with Linked Boundary Conditions},
  journal   = {Studies in Applied Mathematics},
  volume    = {145},
  number    = {3},
  pages     = {357--396},
  year      = {2020},
  doi       = {10.1111/sapm.12322}
}

@article{Bertin2019,
  author    = {Bertin, Karine and El Kolei, Salima and Klutchnikoff, Nicolas},
  title     = {Adaptive Density Estimation on Bounded Domains},
  journal   = {Annals of the Institute of Statistical Mathematics},
  volume    = {71},
  number    = {5},
  pages     = {1081--1109},
  year      = {2019}
}

@inproceedings{Charikar2023,
  author    = {Charikar, Moses and Kapralov, Michael and Nouri, Navid and Woodruff, David},
  title     = {Fast private kernel density estimation via locality sensitive quantization},
  booktitle = {Proceedings of the 40th International Conference on Machine Learning},
  series    = {Proceedings of Machine Learning Research},
  volume    = {202},
  pages     = {3969--3989},
  publisher = {PMLR},
  address   = {Honolulu, HI, USA},
  year      = {2023}
}

@article{Kenig2024,
  author    = {Kenig, Ori and Todros, Koby and Adali, Tulay},
  title     = {Robust {Gaussian} Mixture Modeling},
  journal   = {{IEEE} Transactions on Signal Processing},
  volume    = {72},
  pages     = {3525--3540},
  year      = {2024}
}

@article{Lee2012,
  title   = {{EM} algorithms for multivariate {Gaussian} mixture models with truncated and censored data},
  author  = {Lee, Gyemin and Scott, Clayton},
  journal = {Computational Statistics \& Data Analysis},
  volume = {56},
  number = {9},
  pages = {2816--2829},
  year = {2012},
  issn = {0167-9473},
}

@inproceedings{Bottou2010,
author = {Bottou, Léon},
year = {2010},
month = {09},
pages = {177--186},
title = {Large-Scale Machine Learning with Stochastic Gradient Descent},
isbn = {978-3-7908-2603-6},
booktitle = {Proc. of COMPSTAT},
doi = {10.1007/978-3-7908-2604-3_16}
}

@article{Bottou2018,
  title   = {Optimization Methods for Large-Scale Machine Learning},
  author  = {Bottou, Leon and Curtis, Frank E. and Nocedal, Jorge},
  journal = {SIAM Review},
  volume = {60},
  number = {2},
  pages = {223-311},
  year = {2018},
  doi = {10.1137/16M1080173},
}

@article{Schwartz1978,
  title   = {Estimating the Dimension of a Model},
  author  = {Schwarz, Gideon},
  journal = {The Annals of Statistics},
  volume  = {6},
  year    = {1978},
  number = {2},
  pages = {461--464}
}

@inproceedings{Arthur2007,
author = {Arthur, David and Vassilvitskii, Sergei},
year = {2007},
month = {01},
pages = {1027-1035},
title = {K-Means++: The Advantages of Careful Seeding},
volume = {8},
booktitle = {Proc. of the Annu. ACM-SIAM Symp. on Discrete Algorithms},
doi = {10.1145/1283383.1283494}
}

@inproceedings{Loschilov2017,
  title     = {{SGDR}: Stochastic Gradient Descent with Warm Restarts},
  author    = {Loshchilov, Ilya and Hutter, Frank},
  booktitle = {ICLR},
  year      = {2017}
}

@inproceedings{Liu2022,
  title     = {Unsupervised Anomaly Detection by Robust Density Estimation},
  author    = {Liu, Bo and Tan, Pang-Ning and Zhou, Jiayu},
  booktitle = {AAAI},
  year      = {2022},
  volume    = {36},
  issue     = {4},
  pages     = {4101--4108}
}

@article{Rajwade2008,
  title   = {Probability Density Estimation Using Isocontours and Isosurfaces},
  author  = {Rajwade, Ajit and Banerjee, Arunava and Rangarajan, Anand},
  journal = {IEEE Transactions on Pattern Analysis and Machine Intelligence},
  volume  = {31},
  year    = {2009},
  issue   = {3},
  pages   = {475--491}
}

@book{Gramacki2018,
   title={Nonparametric kernel density estimation and its computational aspects},
  author={Gramacki, Artur},
  volume={37},
  year={2018},
  publisher={Springer},
  address = {Cham, Switzerland}
}

@article{Lin2006,
 ISSN = {10170405, 19968507},
 URL = {http://www.jstor.org/stable/24307552},
 author = {Yi Lin and Ming Yuan},
 journal = {Statistica Sinica},
 number = {2},
 pages = {425--439},
 publisher = {Institute of Statistical Science, Academia Sinica},
 title = {CONVERGENCE RATES OF COMPACTLY SUPPORTED RADIAL BASIS FUNCTION REGULARIZATION},
 volume = {16},
 year = {2006}
}

@article{Wong2002,
title = {Compactly supported radial basis functions for shallow water equations},
journal = {Applied Mathematics and Computation},
volume = {127},
number = {1},
pages = {79--101},
year = {2002},
issn = {0096-3003},
doi = {10.1016/S0096-3003(01)00006-6},
url = {https://www.sciencedirect.com/science/article/pii/S0096300301000066},
author = {S.M. Wong and Y.C. Hon and M.A. Golberg}
}

@article{Morel2018,
title = {Surface reconstruction of incomplete datasets: A novel {Poisson} surface approach based on {CSRBF}},
journal = {Computers \& Graphics},
volume = {74},
pages = {44--55},
year = {2018},
issn = {0097-8493},
doi = {10.1016/j.cag.2018.05.004},
url = {https://www.sciencedirect.com/science/article/pii/S0097849318300633},
author = {Jules Morel and Alexandra Bac and Cédric Véga}
}

@article{Mohammadi2019,
  author  = {Mohammadi, R. and Mahdavi, S. R. and Jaberi, R. and Siavashpour, Z. and Janani, L. and Meigooni, A. S. and Reiazi, R.},
  title   = {Evaluation of deformable image registration algorithm for determination of accumulated dose for brachytherapy of cervical cancer patients},
  journal = {Journal of Contemporary Brachytherapy},
  year    = {2019},
  volume  = {11},
  number  = {5},
  pages   = {469--478},
  doi     = {10.5114/jcb.2019.88762}
}

@article{Buhmann2001,
 ISSN = {00255718, 10886842},
 URL = {http://www.jstor.org/stable/2698936},
 author = {M. D. Buhmann},
 journal = {Mathematics of Computation},
 number = {233},
 pages = {307--318},
 publisher = {American Mathematical Society},
 title = {A New Class of Radial Basis Functions with Compact Support},
 volume = {70},
 year = {2001}
}

@article{Wu1995,
  author  = {Wu, Z.},
  title   = {Compactly supported positive definite radial functions},
  journal = {Advances in Computational Mathematics},
  year    = {1995},
  volume  = {4},
  number  = {1},
  pages   = {283--292},
  doi     = {10.1007/BF03177517}
}

@article{Risser2025,
author = {Risser, Mark D. and Noack, Marcus M. and Luo, Hengrui and Pandolfi, Ronald J.},
title = {Compactly-Supported Nonstationary Kernels for Computing Exact {Gaussian} Processes on Big Data},
journal = {Environmetrics},
volume = {36},
number = {8},
pages = {e70054},
doi = {10.1002/env.70054},
year = {2025}
}

@book{Mclachlan2000,
  author    = {McLachlan, Geoffrey J. and Peel, David},
  title     = {Finite Mixture Models},
  publisher = {John Wiley \& Sons},
  address   = {New York},
  year      = {2000}
}

@book{Scrucca2023,
  author    = {Scrucca, Luca and Fraley, Chris and Murphy, T. Brendan and Raftery, Adrian E.},
  title     = {Model-Based Clustering, Classification, and Density Estimation Using {mclust} in {R}},
  publisher = {Chapman \& Hall/CRC},
  address   = {Boca Raton, FL},
  year      = {2023}
}

@book{Nocedal2006,
  title={Numerical optimization},
  author={Nocedal, Jorge and Wright, Stephen J},
  year={2006},
  publisher={Springer},
  address = {New York}
}

@misc{abalone,
  author = {Nash, Warwick and Sellers, Tracy and Talbot, Simon and Cawthorn, Andrew and Ford, Wes},
  title  = {Abalone},
  year   = {1994},
  howpublished = {UCI Machine Learning Repository},
  doi    = {10.24432/C55C7W}  
}

@misc{glass,
  author = {German, B.},
  title  = {Glass Identification},
  year   = {1987},
  howpublished = {UCI Machine Learning Repository},
  doi    = {10.24432/C5WW2P}  
}

@misc{chickenpox,
  author = {Rozemberczki, Benedek and Scherer, Paul and Kiss, Oliver and Sarkar, Rik and Ferenci, Tamas},
  title  = {{Hungarian} Chickenpox Cases},
  year   = {2021},
  howpublished = {UCI Machine Learning Repository},
  doi    = {10.24432/C5103B} 
}

@article{Johnson1996, 
  author  = {Johnson, Roger W.},
  title   = {Fitting percentage of body fat to simple body measurements},
  journal = {Journal of Statistics Education},
  volume  = {4},
  number  = {1},
  year    = {1996},
  doi     = {10.1080/10691898.1996.11910505} 
}

\end{document}